\documentclass[11pt]{article}
\usepackage{bm}
\usepackage{latexsym,amsfonts}
\usepackage{amssymb,amsmath}
\usepackage{enumitem}
\usepackage{graphics}
\usepackage[demo]{graphicx}
\usepackage{epsfig}
\usepackage{pstricks}
\usepackage[normal]{subfigure}
\usepackage[latin1]{inputenc}
\usepackage[english,francais]{babel}
\usepackage{relsize,exscale}
\usepackage{makeidx}
\usepackage{enumitem}
\usepackage{amsfonts,amssymb,amsmath}
\usepackage{graphicx}
\usepackage{color}
\usepackage{multirow}
\usepackage{mathrsfs}
\usepackage[normalem]{ulem}
\usepackage{cancel}
\usepackage{bbm}
\usepackage{empheq}
\newenvironment{prooff}{{\it Proof :}}{\hfill\rule{2mm}{2mm}\vskip3mm\par}
\newtheorem{theorem}{Theorem}[section]
\newtheorem{lemma}[theorem]{Lemma}
\newtheorem{proposition}[theorem]{Proposition}

\newtheorem{e-definition}[theorem]{Definition\rm}

\usepackage{color}
\definecolor{dred}{rgb}{0.92,0,0}
\definecolor{dgreen}{rgb}{0,0.92,0}
\definecolor{dblue}{rgb}{0,0,0.92}
\definecolor{dyellow}{rgb}{0.95,0.95,0}

\newcommand{\R}{\mathbb{R}}

\def\D{\displaystyle}
\newcommand{\ds}{\displaystyle}
\newcommand{\hs}{\hspace{0.1cm}}
\newcommand{\vs}{\vspace{0.4cm}}
\newcommand{\sa}{\\ [0.2cm]}
\usepackage{multirow}
\graphicspath{
{./Figures/}
{Figures/}
{./}
}
\title{An optimal first-order Taylor-like formula \\ with a minimized remainder}
\author{Jo\"el Chaskalovic \thanks{
D'Alembert, Sorbonne University, Paris, France, (Email:  jch1826@gmail.com)}
\qquad
Franck Assous
\thanks{
Department of Mathematics, Ariel University, Ariel, Israel, (Email: assous@ariel.ac.il)}
\qquad
}
\date{}
\begin{document}
\maketitle
\selectlanguage{english}
\begin{abstract}
\noindent 
In this paper, we derive an optimal first-order Taylor-like formula. In a seminal paper \cite{arXiv_First_Order}, we introduced a new first-order Taylor-like formula that yields a reduced remainder compared to the classical Taylor's formula. Here, we relax the assumption of equally spaced points in our formula. Instead, we consider a sequence of unknown points and a sequence of unknown weights. Then, we solve an optimization problem to determine the best distribution of points and weights that ensures that the remainder is as minimal as possible.
\end{abstract}
\noindent {\em keywords}: Taylor's theorem, Taylor-like formula, Error estimate, Interpolation error, Approximation error, Finite elements.

\section{Introduction}\label{intro}
\noindent Even today, improving the accuracy of approximations remains a challenging problem in numerical analysis. In this context, Taylor's formula plays a crucial role in various domains, especially when one considers error estimates in numerical analysis to assess the accuracy of a numerical approximation method (for example,  see  \cite{RaTho82} for finite element methods).\\

\noindent This challenge becomes even more crucial when comparing the relative accuracy between two given numerical methods. All error estimates share a common structure, whether applied to the finite elements method \cite{BrSc08}, \cite{Ern_Guermond}, numerical approximations of ordinary differential equations \cite{Crouzeix}, or to quadrature formulas used for approximating integrals \cite{Axioms_JCH}.\\

\noindent Let us specify these ideas in this context of numerical integration. Consider, for instance, a composite quadrature rule of order $k$. For a given interval $[a,b]$, let $f$ be a function in $C^{k+1}([a,b])$. The corresponding error of the composite quadrature rule can be expressed as (refer to, e.g., \cite{Atki88}, \cite{BuFa11} or \cite{Crouzeix}), for a non-zero integer $N$:
$$
\D \left|\int_{a}^{b}f(x)dx - \sum_{i=0}^{N}\lambda_i f(x_i)\right| \, \leq \, C_k\,h^{k+1}\,.
$$
In this formula, $h$ denotes the size of the $N+1$ equally spaced panels $[x_{i},x_{i+1}]$, $0 \leq i \leq N$, that discretize the interval $[a,b]$, and  $\lambda_i$ are $N+1$ real numbers. Moreover, $C_k$ is an unknown constant, independent of $h$, but dependent on $f$ and $k$.  This constant is directly linked to the uncertainty associated with the remainder of Taylor's formula.\\

\noindent Usually, to overcome the lack of information associated with the unknown value of the left-hand side  which belongs to the interval $[0, C_k h^{k+1}]$, only the asymptotic convergence rate comparison is considered. This comparison allows us to assess the relative accuracy between two numerical quadratures of order $k_1$ and $k_2, (k_1<k_2)$, when $h$ tends to zero.\sa

\noindent However, comparing the asymptotic convergence rate is no longer useful when evaluating two composite quadratures rules for a fixed value of $h$, as is common in various applications. Therefore, we focus our attention to minimize the constants $C_k$ by improving the estimation of the remainder involved in Taylor's formula.\\

\noindent  In this context, several approaches have been proposed to determine a way to enhance the accuracy of approximation. For example, within the framework of numerical integration, we refer the reader to \cite{Barnett_Dragomir}, \cite{Cerone} or \cite{Dragomir_Sofo}, and references therein. From another point of view, due to the lack of information, heuristic methods were considered, basically based on a probabilistic approach, see for instance \cite{Abdulle}, \cite{AsCh2014}, \cite{Hennig}, \cite{Oates} or \cite{Axioms_JCH}, \cite{CMAM2} and \cite{ChAs20}. This allows to compare different numerical methods, and more precisely finite element, for a given fixed mesh size, \cite{MMA2021}.\\

\noindent In this context, we recently developed a first-order Taylor-like formula in \cite{arXiv_First_Order} and a second-order Taylor-like formula in \cite{JCAM2023}. The goal was to minimize the corresponding remainder by transferring part of the numerical weight of this remainder to the polynomial involved in the Taylor expansion.\\

\noindent In both of these cases, we {\em a priori} introduced a linear combination of $f'$ (and $f''$ in \cite{JCAM2023}) computed at equally spaced points in $[a,b]$, and we determined the corresponding weights in order to minimize the remainder. We demonstrated that the associated upper bound in the error estimate is $2n$ times smaller than the classical one for the first-order Taylor's theorem, and $3/16n^2$ times smaller than the classical second-order Taylors's formula.\\

\noindent In this paper, we relax the assumption of equally spaced points and consider a sequence of unknown points in the interval $[a,b]$, where a given function $f$ needs to be evaluated. Simultaneously, we introduce a sequence of unknown weights to be determined with the goal of minimizing the remainder. Then, we will prove that the remainder of the corresponding first-order expansion is minimized when the points between $a$ and $b$ are equally spaced.\\

\noindent The paper is organized as follows. In Section \ref{B}, we present a new first-order Taylor-like formula built on a sequence of given points $x_k, (k=1,\dots,n-1),$ in $[a,b]$, and given weights $\omega_k, (k=0,\dots,n)$. In Section \ref{Optimal} we derive the main results of this paper which deal with the optimal choice 
of points $x_k$ and weights $\omega_k$ that enable us to minimize the corresponding remainder of the first-order Taylor-like formula. Concluding remarks follow.

\section{A first-order Taylor-like theorem}
\label{B}
\noindent Let us begin by recalling the well-known first order Taylor formula \cite{Taylor}. For two reals $a$ and $b$, $a < b$, we consider a function $f \in \mathcal{C}^{2}([a,b])$. Then, there exists $(m_{2}, M_{2}) \in \mathbb{R}^{2} $ such that, for all $x \in [a,b]$,
$$
m_{2} \leqslant f''(x) \leqslant M_{2}
$$
and we have
\begin{equation}\label{First_Order_Taylor}
f(b) = f(a) + (b-a)f'(a) + (b-a)\epsilon_{a,1}(b),
\end{equation}
with
\[ \lim_{b \to a} \epsilon_{a,1}(b) = 0, \]
and
$$
\frac{(b-a)}{2}m_{2} \leqslant \D \epsilon_{a,1}(b) \leqslant \frac{(b-a)}{2}M_{2}.
$$
In a previous paper \cite{arXiv_First_Order}, we derived a first-order Taylor-like formula with the aim of minimizing the classical remainder $\epsilon_{a,1}(b)$. More precisely, we proved
the following result:
\begin{theorem}\label{theorem_1}
Let $f$ be a real mapping defined on $[a,b]$ which belongs to $\mathcal{C}^{2}([a,b])$, such that: $\forall x \in [a,b], -\infty < m_2 \leqslant f''(x) \leqslant M_2 < +\infty$.\sa
Then, for a given non-zero integer $n$, we have the following first-order expansion:
\begin{equation}\label{Generalized_Taylor}
f(b) = f(a) + (b-a)\left(\frac{f'(a) + f'(b)}{2n} + \frac{1}{n}\sum \limits_{k=1}^{n-1} f'\left(a + k\frac{(b-a)}{n}\right)\right) + (b-a)\epsilon_{a,n+1}(b),
\end{equation}
where
$$
\D |\epsilon_{a,n+1}(b)| \leqslant \frac{(b-a)}{8n}(M_2-m_2).\vs
$$
Moreover, for an {\em a priori} choice of regularly spaced points $\D a + k\frac{(b-a)}{n}$, the remainder $\epsilon_{a,n+1}(b)$ is minimum.
\end{theorem}
In the following of the paper, our goal is to relax the assumption of {\em a priori} equidistant points to determine the optimal set of points $(x_k)_{k=0,n} \in [a,b]$, along with the associated weights $(\omega_k)_{k=0,n}$. This determination  will enable us to minimize the remainder $\epsilon_{a,n+1}(b)$ in (\ref{Generalized_Taylor}).\sa
To derive the main result below, we first introduce the function $\phi$ defined by
\begin{eqnarray}
\D \phi: & [0,1] & \!\!\longrightarrow \hs \R \nonumber \\[0.2cm]
\D & t  &  \!\! \longmapsto \hs f'(a + t(b-a))\,, \label{phi}
\end{eqnarray}
that satisfies $\phi(0) = f'(a)$ and $\phi(1) = f'(b)$.
Moreover, we proved in \cite{arXiv_First_Order} that the remainder $\epsilon_{a,1}(b)$ introduced in (\ref{First_Order_Taylor}) satisfies the following result:
\begin{proposition}
The remainder $\epsilon_{a,1}(b)$ in formula (\ref{First_Order_Taylor}) can be expressed as
\begin{equation}\label{Epsilon1}
\epsilon_{a,1}(b) = \int_{0}^{1}{(1-t)\phi'(t)dt}\,.
\end{equation}
\end{proposition}
For a given integer $n \in \mathbb{N}^{*}$, we consider the set of points $(x_k)_{k=0,n}$ in the interval $[a,b]$, where $x_0=a$ and $x_n=b$, on the one hand, and a set of real
weights $(\omega_k)_{k=0,n}$, on the other hand.\sa
Then, we define the quantity $\epsilon^{*}_{n}(a,b)$ by the formula
\begin{equation}\label{Def_epsilon_n+1}
f(b) = f(a) + (b-a)\left(\sum \limits_{k=0}^{n}\omega_{k}f'(x_k)\right) + (b-a)\epsilon^{*}_{n}(a,b),
\end{equation}
where the two sequences $(x_k)_{k=1,n-1}$ and $(\omega_k)_{k=0,n}$  are to be determined to minimize the remainder $\epsilon^{*}_{n}(a,b)$.\sa
For the upcoming result, we need to introduce some notations. We denote by $t_k, (k=0,\dots, n),$ a sequence of real numbers that allows us  to represent the points $x_k$ in $[a,b]$ as a barycentric combination of $a$ and $b$, that is:
\begin{equation}\label{barycentric}
\D x_k = a +t_k(b-a),\, (0\leq t_k\leq 1).
\end{equation}
We also introduce $S_k$ as the partial sum of the weights $\omega_j$, $0 \leq j \leq k$, 
\begin{equation}
\label{defSk}
\D S_k = \sum \limits_{j=0}^{k}\omega_j.
\end{equation}
Now, we can prove the following result:
\begin{theorem}\label{theorem_2}
Let $f$ be a real mapping defined on $[a,b]$ which belongs to $\mathcal{C}^{2}([a,b])$. Then, the remainder $\epsilon^{*}_{n}(a,b)$ defined by (\ref{Def_epsilon_n+1}) satisfies:
\begin{equation}\label{F8_0}
\D\frac{(b-a)}{2}\sum_{k=0}^{n-1}\bigg[m_2(S_k-t_k)^2 - M_2(S_k - t_{k+1})^2\bigg] \leq \epsilon^{*}_{n}(a,b) \leq \frac{(b-a)}{2}\sum_{k=0}^{n-1}
\bigg[M_2(S_k-t_k)^2 - m_2(S_k - t_{k+1})^2\bigg].
\end{equation}
\end{theorem}
\begin{prooff}
Using first the function $\phi$ defined in (\ref{phi}), formulas (\ref{First_Order_Taylor}) and (\ref{Def_epsilon_n+1}) lead to
\[ \D\frac{f(b) - f(a)}{b - a} = \phi(0) + \epsilon_{a,1}(b) = \sum \limits_{k=0}^{n} \omega_{k}\phi(t_k) + \epsilon^{*}_{n}(a,b),\]
which can be re-written as:
\begin{eqnarray}
\epsilon^{*}_{n}(a,b) & = &\phi(0) + \epsilon_{a,1}(b) - \sum \limits_{k=0}^{n} \omega_{k}\phi(t_k), \nonumber \\[0.2cm]
& = & \phi(0) + \int_{0}^{1}{(1-t)\phi'(t)dt} - \sum \limits_{k=0}^{n} \omega_{k}\phi(t_k),\nonumber\\[0.2cm]
& = & \phi(1) - \int_{0}^{1}{t\phi'(t)dt} - \sum \limits_{k=0}^{n} \omega_{k}\phi(t_k), \nonumber\\[0.2cm]
& = &\phi(1) - \int_{0}^{1}{t\phi'(t)dt} + \sum \limits_{k=0}^{n} \omega_{k}\left(\phi(1)-\phi(t_k)\right) - \sum \limits_{k=0}^{n} \omega_{k}\phi(1), \nonumber \\[0.2cm]
& = &\left(1 - \sum \limits_{k=0}^{n}\omega_{k} \right)\phi(1) - \int_{0}^{1}{t\phi'(t)dt} + \sum \limits_{k=0}^{n} \omega_{k}\left(\phi(1)-\phi(t_k)\right).\label{Equation espilon_n+1_1}
\end{eqnarray}
Let us assume that weights $(\omega_k)_{k=0,n}$ fulfill the condition
\begin{equation}\label{poids_norme}
\D \sum \limits_{k=0}^{n} \omega_{k} = 1\,,
\end{equation}
equation (\ref{Equation espilon_n+1_1}) can be expressed as
\begin{eqnarray}
\epsilon^{*}_{n}(a,b) & = & - \int_{0}^{1}{t\phi'(t)dt} + \sum\limits_{k=0}^{n} \omega_{k}\left(\phi(1)-\phi(t_k)\right),\nonumber \\[0.2cm]
 & = & - \int_{0}^{1}{t\phi'(t)dt} + \sum \limits_{k=0}^{n}\omega_{k}\int_{t_k}^{1}{\phi'(t)dt}, \nonumber \\[0.2cm]
 & = & - \sum \limits_{k=0}^{n-1}\int_{t_k}^{t_{k+1}}{t\phi'(t)dt} + \sum \limits_{k=0}^{n}\omega_{k}\int_{t_k}^{1}{\phi'(t)dt}.\label{Equation espilon_n+1_2}
\end{eqnarray}
Considering the second term of (\ref{Equation espilon_n+1_2}), it can be transform as follows:
\begin{eqnarray}
\hspace{-0.7cm}\D \sum \limits_{k=0}^{n}\omega_{k}\int_{t_k}^{1}\phi'(t)dt &\hspace{-0.2cm} = & \hspace{-0.1cm}\int_{t_0}^{1}\omega_{0}\phi'(t)dt + \int_{t_1}^{1}\omega_{1}\phi'(t)dt + \dots
+ \int_{t_{n-1}}^{1}\omega_{n-1}\phi'(t)dt, \nonumber \\[0.2cm]
\hspace{-0.7cm}&\hspace{-0.2cm} = & \hspace{-0.1cm}\int_{t_0}^{t_1}\!\!\omega_{0}\phi'(t)dt + \int_{t_1}^{t_2}\!(\omega_{0}+\omega_{1})\phi'(t)dt + \dots + \int_{t_{n-1}}^{1}\!\!\!\!(\omega_{0}+\dots+\omega_{n-1})\phi'(t)dt. \nonumber \\[0.2cm]
\hspace{-0.7cm}&\hspace{-0.2cm} = & \hspace{-0.1cm}\sum \limits_{k=0}^{n-1}\int_{t_k}^{t_{k+1}}\!S_{k}\,\phi'(t)dt, \label{New_F1}
\end{eqnarray}
where $S_k$ is defined by (\ref{defSk}).\sa
Then, using (\ref{New_F1}) in (\ref{Equation espilon_n+1_2}) enables to write $\epsilon^{*}_{n}(a,b)$ as
\begin{equation}\label{F1}
\D\epsilon^{*}_{n}(a,b) = \sum \limits_{k=0}^{n-1} \int_{t_k}^{t_{k+1}}{(S_{k} - t)\phi'(t)dt}.
\end{equation}
In the estimations below, we will use that
$$
\forall\, t \in [0,1], \phi'(t) = (b-a) f''(a + t(b-a)), \hs \mbox{ and } \forall \,x \in [a,b],  m_2 \leqslant f''(x) \leqslant M_2.
$$
Next, to derive a double inequality on $\epsilon^{*}_{n}(a,b)$, we consider the three following cases, depending on the location of $S_k$ related to the interval $[t_k, t_{k+1}]$.\\
\begin{enumerate}
\item If $t_k \leq S_k \leq t_{k+1}$, the integral in (\ref{F1}) can be decomposed as follows:
\begin{equation}
\int_{t_k}^{t_{k+1}}{(S_{k} - t)\phi'(t)dt} = \int_{t_k}^{S_{k}}(S_{k} - t)\phi'(t) dt + \int_{S_{k}}^{t_{k+1}}(S_{k} - t)\phi'(t) dt.
\end{equation}
Now, using that $(S_{k}-t)$ is positive on $[t_k, S_{k}]$, and negative on $[S_{k},t_{k+1}]$, we can write
$$
\left\{
\begin{array}{l}
(b-a)m_2 \ds \int_{t_k}^{S_{k}}(S_{k} - t) dt \leq \int_{t_k}^{S_{k}}(S_{k} - t)\phi'(t) dt \leq (b-a)M_2\int_{t_k}^{S_{k}}(S_{k} - t) dt,\\
\\
(b-a)M_2 \ds\int_{S_{k}}^{t_{k+1}}(S_{k} - t) dt \leq \int_{S_{k}}^{t_{k+1}}(S_{k} - t)\phi'(t) dt \leq (b-a)m_2\int_{S_{k}}^{t_{k+1}}(S_{k} - t) dt.
\end{array}
\right.
$$
Summing up these two relations, we obtain first that
$$
\int_{t_k}^{t_{k+1}}(S_{k}-t)\phi'(t) dt \leq (b-a)M_2\int_{t_k}^{S_{k}}(S_{k}-t) dt + (b-a)m_2\int_{S_{k}}^{t_{k+1}}(S_{k}-t) dt,
$$
and also
$$
\int_{t_k}^{t_{k+1}}(S_{k}-t)\phi'(t) dt \geq (b-a)m_2\int_{t_k}^{S_{k}}(S_{k}-t) dt + (b-a)M_2\int_{S_{k}}^{t_{k+1}}(S_{k}-t) dt\,.
$$
By simply computing the integrals involved in these relations, we obtain, in that case,  the estimate
\begin{equation}
\label{F67}
\frac{(b-a)}{2}\bigg( \! m_2(S_k-t_k)^2 - M_2(S_k - t_{k+1})^2 \! \bigg)
 \!\!\leq \!\!
 \int_{t_k}^{t_{k+1}} \!\!(S_{k}-t)\phi'(t) dt 
 \! \leq  \!\!
 \frac{(b-a)}{2}\bigg( \! M_2(S_k-t_k)^2 - m_2(S_k - t_{k+1})^2 \! \bigg).
\end{equation}
\item If $S_k \geq t_{k+1}$ then $S_k - t \geq 0$ for all $t \in [t_k,t_{k+1}]$, and we have:
$$
\D (b-a)m_2\int_{t_k}^{t_{k+1}}(S_k - t) dt \leq \int_{t_k}^{t_{k+1}}(S_{k}-t)\phi'(t) dt \leq (b-a)M_2\int_{t_k}^{t_{k+1}}(S_k - t) dt\,.
$$
This yields
$$
\D \frac{(b-a)m_2}{2}\bigg((S_k-t_k)^2 - (S_k - t_{k+1})^2\bigg) \leq \int_{t_k}^{t_{k+1}}(S_{k}-t)\phi'(t) dt \leq \frac{(b-a)M_2}{2}\bigg((S_k-t_k)^2 - (S_k - t_{k+1})^2\bigg),
$$
which also leads to (\ref{F67}), by simply using that $m_2 \leq M_2$.\\
\item If $S_k \leq t_{k}$ then $S_k - t \leq 0$ for all $t \in [t_k,t_{k+1}]$, and in the same way as above, we get
$$
\D \frac{(b-a)M_2}{2}\bigg((S_k-t_k)^2 - (S_k - t_{k+1})^2\bigg) \leq \int_{t_k}^{t_{k+1}}(S_{k}-t)\phi'(t) dt \leq \frac{(b-a)m_2}{2}\bigg((S_k-t_k)^2 - (S_k - t_{k+1})^2\bigg),
$$
that also gives estimate (\ref{F67}).
\end{enumerate}
Hence, in all cases, we arrive at the same estimate (\ref{F67}). Finally, by summing over (\ref{F67}) over all values of $k$ from $0$ to $n-1$ in  (\ref{F67}), we get inequalities (\ref{F8_0}) for the remainder $\epsilon^{*}_{n}(a,b)$, which concludes the proof.
\end{prooff}
With the aim of minimizing $\epsilon^{*}_{n}(a,b)$, we introduce the function $\chi$ defined by
\begin{equation}\label{chi}
\D\chi = \frac{(b-a)(M_2-m_2)}{2}\sum_{k=0}^{n-1}\bigg[\bigg(\sum \limits_{j=0}^{k}\omega_j -t_k \bigg)^{\!\!2} + \bigg(\sum \limits_{j=0}^{k}\omega_j - t_{k+1}\bigg)^{\!\!2}\,\bigg],
\end{equation}
which represents to the difference between the right-hand side and the left-hand side in (\ref{F8_0}). \sa
As a consequence, in the next section, we will minimize function $\chi$ which depends on $2n-1$ variables, namely $(t_1,\dots,t_{n-1},\omega_0,\dots,\omega_{n-1})$. This is because the boundary points are known  ($x_0=a$ corresponding to $t_0=0$ and $x_{n}=b$ corresponding to $t_n=1$), on the one hand, and due to relation (\ref{poids_norme}) which will determine $\omega_n$, on the other hand.
\section{The optimal first order Taylor-like formula}\label{Optimal}
\noindent We begin this section by deriving a lemma that provides a necessary condition for the function $\chi$ to have an extremum at the point $(t_1,\dots,t_{n-1},\omega_0,\dots,\omega_{n-1})$.
\begin{lemma}\label{Min1}
Let $(t_1,\dots,t_{n-1},\omega_0,\dots,\omega_{n-1})$ be an extremum of function $\chi$. \sa
Then, we have:
\begin{empheq}[left=\mbox{}\hs\empheqlbrace]{alignat=2}
\D \hs S_k & = \hs\frac{1}{2}(t_k+t_{k+1}), & \hs (k=0,\dots,n-1), \label{Extrem_1}\\[0.2cm]
\D \hs t_k & = \hs S_{k-1} + \frac{\omega_k}{2}, & (k=1,\dots,n-1). \label{Extrem_2}
\end{empheq}
\end{lemma}
\begin{prooff}
The necessary conditions that guarantee that $(t_1,\dots,t_{n-1},\omega_0,\dots,\omega_{n-1})$ is an extremum is written as (see for instance \cite{extremum})
$$
\D \forall k=0,\dots,n-1: \frac{\partial\chi}{\partial \omega_k} = 0 \hs \mbox{ and } \hs \forall k=1,\dots, n-1: \frac{\partial\chi}{\partial t_k} = 0.
$$
Regarding first the dependence of the function $\chi$ on the variables $\omega_k$, $(k=0,\dots,n-1)$, the conditions $\D \frac{\partial\chi}{\partial \omega_k} = 0$ are expressed as
$$
\sum_{m=k}^{n-1}\bigg[2S_m-t_m-t_{m+1}\bigg]=0\,,
$$
or equivalently
\begin{equation}
\label{XX}
\ds \sum_{m=k}^{n-1}S_m = \frac{1}{2}\sum_{m=k}^{n-1}(t_m+t_{m+1})\,.
 \end{equation}
Since this system of equations is triangular, it can be easily solved. Writing two consecutive equations for a given $k \in \{0,\dots,n-2 \}$ leads to 
\begin{eqnarray*}
\D S_k & + & S_{k+1} +  \dots + S_{n-1} \, = \, \frac{1}{2}\big[(t_k + t_{k+1}) + (t_{k+1} + t_{k+2}) + \dots + (t_{n-1} + t_n)\big], \\[0.2cm]
\D  &  & S_{k+1} +  \dots + S_{n-1} \, = \, \frac{1}{2}\big[(t_{k+1} + t_{k+2}) + \dots + (t_{n-1} + t_n)\big]\,, 
\end{eqnarray*}
which readily gives, by difference,
$$
\D S_k = \frac{1}{2}(t_k + t_{k+1})\,,
$$
the case $k=n-1$ corresponding directly to the last equation of the system (\ref{XX}).\\

\noindent Now, to study the dependence of the function $\chi$ on the variables $t_k$, $\forall k=1,\dots,n-1$, we  expand formula (\ref{chi}), using that  $t_0=0$ and $t_n=1$. We obtain that
\begin{eqnarray*}
\D\frac{2}{(b-a)(M_2-m_2)}\, \chi & = & \sum_{k=0}^{n-1}\bigg[\bigg(\omega_0 + \dots +\omega_k - t_k \bigg)^{\!\!2} + \bigg(\omega_0 + 
\dots +\omega_k - t_{k+1}\bigg)^{\!\!2}\,\bigg] \nonumber \\[0.2cm]
\D  & = & (\omega_0 - 0)^2 + (\omega_0 - t_1)^2   \nonumber \\[0.2cm]
\D &  &+ \,(\omega_0 + \omega_1 - t_1)^2 + (\omega_0 + \omega_1 - t_2)^2  \nonumber \\[0.2cm]
\D &  & +\, \dots \nonumber \\[0.2cm]
\D &  & + \,(\omega_0 + \dots + \omega_{k} - t_{k-1})^2 + (\omega_0 + \dots + \omega_{k} - t_{k})^2   \nonumber \\[0.2cm].
\D &  & +\, \dots \nonumber \\[0.2cm]
\D &  & + \,(\omega_0 + \dots + \omega_{n-1} - t_{n-1})^2 + (\omega_0 + \dots + \omega_{n-1} - 1)^2\,.
\end{eqnarray*}

\noindent So, by taking the derivative of the function $\chi$ with respect to $t_k$, we obtain, for each $k=1,\dots,n-1$: 
\begin{eqnarray*}
\D \frac{\partial\chi}{\partial t_k} = 0 & \Leftrightarrow & 2(\omega_0 + \dots + \omega_{k-1} - t_k) + 2 (\omega_0 + \dots + \omega_k - t_k) = 0\,.
\end{eqnarray*}
This can be expressed as
$$
S_{k-1} + S_k - 2t_k =0
$$
that is, using the definition (\ref{defSk}) of $S_k$
$$
2t_k = 2S_{k-1} +  \omega_k\,,
$$
which corresponds to (\ref{Extrem_2}).
\end{prooff}
From Lemma \ref{Min1}, we can now state the main result of this paper:
\begin{theorem}
Let $f$ be a real function defined on $[a,b]$ that belongs to $\mathcal{C}^{2}([a,b])$. Then, the optimal unknown weights $(\omega_m)_{m=0,n}$ together with the optimal
set of points $(x_m)_{m=1,n-1}$ determined by the sequence of real numbers $(t_m)_{m=1,n-1}$ that minimizes the remainder $\epsilon^{*}_{n}(a,b)$ defined by
(\ref{Def_epsilon_n+1}) are given by:
\begin{equation}\label{poids_thm}
\D \omega_{0} = \omega_{n} = \frac{1}{2n} \,\mbox{ and } \,\omega_{k} = \frac{1}{n}, \forall k=1,\dots,n-1,
\end{equation}
\begin{equation}\label{tk_thm}
\D t_k = \frac{k}{n} \,\mbox{ and }\, x_k = a  + \frac{k}{n}(b-a), \,\forall k=1,\dots,n-1.
\end{equation}
As a result, the corresponding optimal first-order Taylor-like formula is given by the following expression:
\begin{equation}\label{Generalized_Taylor_2}
\D f(b) = f(a) + (b-a)\left(\frac{f'(a) + f'(b)}{2n} + \frac{1}{n}\sum \limits_{k=1}^{n-1} f'\left(a + k\frac{(b-a)}{n}\right)\right) + (b-a)\epsilon^{*}_{n}(a,b),
\end{equation}
with
\begin{equation}\label{majoration_epsilon_n+1_2}
\D |\epsilon^{*}_{n}(a,b)| \leqslant \frac{(b-a)}{8n}(M_2-m_2).
\end{equation}
\end{theorem}
\begin{prooff}
\begin{itemize}
\item We begin by proving that $\omega_k$ is constant for all values of $k$ belonging to $\{1,\dots,n-1\}$ \sa
From (\ref{Extrem_1}) and (\ref{Extrem_2}), we have
$$
\D 2 S_k = t_k + t_{k+1} = S_{k-1} + \frac{\omega_k}{2} + S_{k} + \frac{\omega_{k+1}}{2}, \, \forall k=1,\dots, n-2,
$$
that yields
$$
\D S_{k} - S_{k-1} := \omega_k = \frac{\omega_{k}+\omega_{k+1}}{2}.
$$
Then,
\begin{equation}\label{Omegak_Const}
\D\omega_{k+1}=\omega_{k}, \, \forall k=1,\dots,n-2,
\end{equation}
that corresponds to
$$
\D \omega_1 = \dots = \omega_{n-1}.
$$
\item We will now establish the relation between $\omega_0$ and $\omega_k$, for $k=1,\dots,n-1$.\sa
Firstly, let us write (\ref{Extrem_1}) for $k=0$ and (\ref{Extrem_2}) for $k=1$. Using that $t_0 = 0$, we obtain that
\begin{eqnarray*}
\D 2S_0 & = & t_0 + t_1 = t_1,  \\ [0.2cm]
\D t_1 & = & S_0 + \frac{\omega_1}{2} = \omega_0 +\frac{\omega_1}{2},
\end{eqnarray*}
from which we get 
$$
\D \omega_1= 2\omega_0\,.
$$
This allows us to conclude, using (\ref{Omegak_Const}), that
\begin{equation}\label{omega_Const}
\D \omega_1 = \omega_2 = \dots = \omega_{n-1} = 2 \omega_0.
\end{equation}
\item Let us compute now the value of $\omega_0$.\sa
To this end, we write (\ref{Extrem_1}) and (\ref{Extrem_2}) for $k=n-1$. Given that $t_n=1$ and using (\ref{omega_Const}), this yields 
\begin{eqnarray*}
\D 2S_{n-1} & = & t_{n-1} + t_n = t_{n-1} +1, \\[0.2cm]
\D & = & S_{n-2} + \frac{\omega_{n-1}}{2} + 1.
\end{eqnarray*}
Then, substituting the expressions of $S_{n-1}$ and $S_{n-2}$, we get with (\ref{omega_Const}) that
$$
\D 2\omega_0 + 4(n-1)\omega_0 = \omega_0 + 2(n-2)\omega_0 + \frac{\omega_{n-1}}{2} + 1,
$$
that leads to
$$
\D 2(2n+1)\omega_0 - 2 = \omega_{n-1} = 2\omega_0
$$
that is
\begin{equation}\label{omega0_Value}
\D \omega_0 = \frac{1}{2n}.
\end{equation}
\item It remains now to compute the value of $\omega_n$. 
Using (\ref{omega_Const}) and (\ref{omega0_Value}) gives directly that
\begin{equation}\label{omega_Const_Val}
\D \omega_1 = \omega_2 = \dots = \omega_{n-1} = \frac{1}{n}.
\end{equation}
Finally, the value of $\omega_n$ is obtained from relation (\ref{poids_norme}), that is
$$
\D \omega_n = 1 - \omega_0 - \sum_{k=1}^{n-1}\omega_k = 1 - \frac{1}{2n} - \frac{n-1}{n} = \frac{1}{2n}.
$$
\item Let us now consider the values of the $t_k$, for $k=1\dots,n-1$.
\sa
Using the expressions (\ref{Extrem_2}) of the optimal $t_k,$ together with the expressions (\ref{omega0_Value}) and (\ref{omega_Const_Val}) of the $\omega_k$, we obtain that
\begin{eqnarray*}
\D t_k & = & \omega_0 + \sum_{j=1}^{k-1}\omega_j + \frac{\omega_k}{2}, (k=1,\dots,n-1)\\[0.2cm]
\D & = & \frac{1}{2n} + \frac{k-1}{n} + \frac{1}{2n} = \frac{k}{n},
\end{eqnarray*}
that yields with (\ref{barycentric}), the following expressions of the optimal points $x_k$:
\begin{equation}
\D x_k = a + \frac{k}{n}(b-a), (k=1,\dots,n-1).
\end{equation}
\item We conclude the proof of this theorem by determining the optimal lower and upper bounds of the remainder $\epsilon^{*}_{n}(a,b)$ given in (\ref{F8_0}).
\begin{enumerate}
\item We first evaluate the quantity $\D\sum_{k=0}^{n-1}(S_k-t_k)^2$:\sa
Using the expression (\ref{poids_thm}) of the $\omega_j$, we obtain that 
$$
\D S_k = \sum_{j=0}^{k}\omega_j = \omega_0 + \sum_{j=1}^{k}\omega_j = \frac{1}{2n} +\frac{k}{n} = \frac{2k+1}{2n}, (k=0,\dots,n-1)\,,
$$
so that, together with the expression (\ref{tk_thm}) of the $t_k$, we have
\begin{equation}\label{Sum1}
\D\sum_{k=0}^{n-1}(S_k-t_k)^2 = \sum_{k=0}^{n-1}\bigg(\frac{2k+1}{2n}-\frac{k}{n}\bigg)^2 = \frac{1}{4n}.
\end{equation}
\item Similarly, we evaluate the quantity $\D\sum_{k=0}^{n-1}(S_k-t_{k+1})^2$. The same arguments yield
\begin{eqnarray}
\D\sum_{k=0}^{n-1}(S_k-t_{k+1})^2  & = & \sum_{k=0}^{n-2}(S_k-t_{k+1})^2 + (S_{n-1}-t_n)^2 \nonumber\\[0.2cm]
\D & = & \sum_{k=0}^{n-2}\bigg(\frac{2k+1}{2n}-\frac{k+1}{n}\bigg)^2 +\bigg(\frac{2n-1}{2n}-1\bigg)^2 \nonumber\\[0.2cm]
\D & = & \sum_{k=0}^{n-2}\bigg(\frac{-1}{2n}\bigg)^2 + \frac{1}{4n^2} = \frac{1}{4n}. \label{Sum2}
\end{eqnarray}
\end{enumerate}
Therefore, combining (\ref{Sum1}) and (\ref{Sum2}), we obtain optimal lower bound and upper bounds in (\ref{F8_0}) for the remainder $\epsilon^{*}_{n}(a,b)$ as follows:
$$
\D\frac{(b-a)}{8n}(m_2- M_2) \leq \epsilon^{*}_{n}(a,b) \leq \frac{(b-a)}{8n}(M_2 - m_2),
$$
which exactly corresponds to the result derived in \cite{arXiv_First_Order}.
\end{itemize}
\end{prooff}
\section{Conclusion}
\label{Conclusion}
\noindent In this paper, we have derived a new first-order Taylor-like formula constructed as a linear combination of the first derivatives of a given function, evaluated at specified points
$x_k, (k=1,\dots,n-1),$ within an interval $[a,b]$. These points are weighted by real numbers $\omega_k, (k=0,\dots,n)$. Unlike the approach in \cite{arXiv_First_Order}, the positions of these points and their corresponding weights were not fixed {\em a priori}. In particular, we relax the assumption of equally spaced points.\\

\noindent Our aim was to determine the optimal positions of these points together with the weights for obtaining the ``best formula'', in the sense that the corresponding remainder is as small as possible. To achieve this, we have established an initial result that provides upper and lower bounds for the remainder.\\

\noindent Then, we proved the existence of an optimal set of points and weights that minimize the remainder in the first-order Taylor-like formula. This result corresponds to the formula we derived in \cite{arXiv_First_Order}, where we explicitely set \emph{a priori} the values of the points $x_k, (k=1,\dots,n-1),$ as uniformly distributed within the interval $[a,b]$.\\

\noindent So, the consequences derived in \cite{arXiv_First_Order}, mainly related to applications in error approximations, can be considered as optimal. Mainly, it treats on $P_1$-interpolation error estimate, the corrected trapezoidal quadrature formula, and finite element error approximations. For example, using the corrected trapezoidal quadrature formula, we have obtained an upper bound which is two times smaller than the errors obtained by using the classical trapezoidal quadrature formula. It highlights the importance and impact of the new Taylor-like formula (\ref{Generalized_Taylor_2})-(\ref{majoration_epsilon_n+1_2}) in assessing the accuracy of a given numerical approximation method.\\

\noindent This research can be extended to the case of dimension strictly greater than one. This extension requires a Taylor-like formula we have already derived in \cite{ChAs2023}. Additionally, we could expand this work to a second-order Taylor-like formula, as proposed in \cite{JCAM2023}. Both of these extensions will be explored to assess their impact on error estimates in the context of applications in numerical analysis.\\

\noindent \textbf{\underline{Homages}:} The authors warmly dedicate this research to pay homage to the memory of Professors Andr\'e Avez and G\'erard Tronel, who largely promote the passion for research and teaching in mathematics of their students.
\end{document}